\documentclass[article]{IEEEtran}

\usepackage{cite}
\usepackage{color}
\usepackage[pdftex]{graphicx}
\graphicspath{{fig/}{jpeg/}}
\usepackage[cmex10]{amsmath}
\usepackage{amssymb}
\usepackage{algorithm}
\usepackage{algorithmic}
\usepackage{amsmath,amssymb,lipsum}
\usepackage{comment}
\usepackage{enumerate}
\usepackage{graphicx}
\usepackage[font=small]{caption}
\usepackage{subcaption}
\usepackage{setspace}
\usepackage{caption}

\input{mysymbol.sty}
\usepackage{needspace}

\usepackage{theorem}
\newtheorem{thm}{Theorem}

\newtheorem{proposition}{Proposition}
\theoremstyle{definition}
\newtheorem{assumption}{Assumption}
\newtheorem{corollary}{Corollary}

\title{Decentralized Quadratically Approximated \\ Alternating Direction Method of Multipliers}

\author{Aryan Mokhtari$^{ \dagger}$ \qquad Wei Shi$^{\star}$ \qquad Qing Ling$^{\star}$ \qquad Alejandro Ribeiro$^{\dagger}$
    \thanks{Supported by NSF CAREER CCF-0952867 and ONR N00014-12-1-0997.}
    
$^{\dagger}$Department of Electrical and Systems Engineering, University of Pennsylvania \\ 
    $^{\star}$Department of Automation, University of Science and Technology of China}

\begin{document}

\maketitle
\thispagestyle{empty}
\pagestyle{empty}

\begin{abstract}

This paper considers an optimization problem that components of the objective function are available at different nodes of a network and nodes are allowed to only exchange information with their neighbors. The decentralized alternating method of multipliers (DADMM) is a well-established iterative method for solving this category of problems; however, implementation of DADMM requires solving an optimization subproblem at each iteration for each node. This procedure is often computationally costly for the nodes. We introduce a decentralized quadratic approximation of ADMM (DQM) that reduces computational complexity of DADMM by minimizing a quadratic approximation of the objective function. Notwithstanding that DQM successively minimizes approximations of the cost, it converges to the optimal arguments at a linear rate which is identical to the convergence rate of DADMM. Further, we show that as time passes the coefficient of linear convergence for DQM approaches the one for DADMM. Numerical results demonstrate the effectiveness of DQM.

%We study the problem of minimizing a sum of convex objective functions where the components of the objective are available at different nodes of a network. The use of distributed subgradient or gradient methods is widespread but they often suffer from slow convergence. In this paper we propose a decentralized network Newton (NN) method that achieves faster convergence. Similar to distributed gradient methods, NN solves a surrogate problem that is an approximate version of the original. The surrogate problem can be efficiently solved using an inexact Newton method that uses an approximation to the Newton step that can be computed in a decentralized manner. Theoretical analysis of NN shows a superlinear rate of convergence to a neighborhood of the optimal solution and demonstrates a tradeoff between accuracy of solutions and speed of convergence. Advantages with respect to distributed gradient descent are further studied through numerical analyses.

\end{abstract} 


%!TEX root = root.tex

\section{Introduction}\label{sec_Introduction}

Decentralized algorithms are designed to solve the problem of minimizing a global cost function over a set of nodes. Agents (nodes) only have access to their local cost functions and try to minimize the global cost cooperatively only by exchanging information with their neighbors. To be more precise, consider a variable $\tbx\in\reals^p$ and a connected network containing $n$ agents each of which has access to a strongly convex local cost function $f_i:\reals^p\to\reals$. The agents cooperate to solve %the optimization problem
\begin{equation}\label{original_optimization_problem1}
  \tbx^*
      \  =\ \argmin_{\tbx}\sum_{i=1}^{n} f_i(\tbx).
\end{equation}
Problems of this form arise in decentralized control systems \cite{Bullo2009,Cao2013-TII,LopesEtal8}, wireless communication networks \cite{Ribeiro10,Ribeiro12}, sensor networks \cite{Schizas2008-1,KhanEtal10,cRabbatNowak04}, and large scale machine learning systems \cite{bekkerman2011scaling,Tsianos2012-allerton-consensus,Cevher2014}. %In this paper we assume that the local costs $f_i$ are twice differentiable and strongly convex. 

There are different algorithms to solve \eqref{original_optimization_problem1} in a distributed manner \cite{Nedic2009,Jakovetic2014-1,YuanQing,wei2012distributed,terelius2011decentralized,Shi2014,NN-part1,NN-part2,Duchi2012,cTsianosEtal12,luo1993convergence,Schizas2008-1,jakovetic2015linear,BoydEtalADMM11,Shi2014-ADMM}. Decentralized implementations of the alternating direction method of multipliers (DADMM) are well-known for solving \eqref{original_optimization_problem1}, with a fast linear convergence rate \cite{Schizas2008-1,BoydEtalADMM11,Shi2014-ADMM,iutzeler2013explicit}. On the other hand, DADMM steps are computationally costly, since each node has to minimize a convex problem at each iteration. This issue is addressed by the Decentralized Linearized ADMM (DLM) algorithm that in lieu of minimizing the exact primal convex problem, minimizes a first-order linearized version of the primal objective function \cite{cQingRibeiroADMM14,ling2014dlm}. DLM reduces the computational complexity of DADMM, however, convergence rate of DLM is slow when the primal objective function is ill-conditioned, since DLM operates on first-order information. Moreover, the coefficient of linear convergence for DLM is strictly smaller than the one for DADMM {\cite{DQMjournal}}. These drawbacks can be resolved by incorporating second-order information of the primal function.

 In this paper we propose a decentralized quadratic approximation of ADMM (DQM) that successively minimizes a quadratic approximation of the primal objective function. Therefore, DQM enjoys the low computational cost of DLM, while its convergence is exact and linear same as DADMM. This approximation incorporates the second-order information of the primal objective function which leads to a provably fast linear convergence even in the case that the cost function is ill-conditioned. Further, we show that as time progresses the coefficient of linear convergence for DQM approaches the one for DADMM.

We begin the paper with describing the DADMM steps and explaining the idea of the DLM algorithm (Section \ref{sec:problem}). Then, we introduce the DQM method which is different from DADMM and DLM in minimizing a quadratic approximation of the primal function (Section \ref{sec:DQM}). We analyze convergence properties of DQM (Section \ref{sec:Analysis}) and evaluate its performance in solving a logistic regression problem (Section \ref{sec:simulations}). Proofs of the results in this paper and comprehensive numerical results are provided in \cite{DQMjournal}.

%!TEX root = root.tex

%%%%%%%%%%%%%%%%%%%%%%%%%%%%%%%%%%%%%%%%%%%%%%%%%%%%%%%%%%%%%%%
%%%   S   E   C   T   I   O   N   %%%%%%%%%%%%%%%%%%%%%%%%%%%%%
%%%%%%%%%%%%%%%%%%%%%%%%%%%%%%%%%%%%%%%%%%%%%%%%%%%%%%%%%%%%%%%
%
\section{DADMM: Decentralized alternating direction method of multipliers} \label{sec:problem}

Consider a connected network with $n$ nodes and $m$ edges where the set of nodes is $\ccalV=\{1, \dots, n\}$ and the set of ordered edges $\ccalE$ contains pairs of nodes $(i,j)$. Node $i$ can communicate with node $j$ if and only if the pair $(i,j)\in \ccalE$ and we further assume that the network is symmetric so that $(i,j)\in \ccalE$ implies $(j,i)\in \ccalE$. We define the neighborhood of node $i$ as the set $\mathcal{N}_i=\{j\mid (i,j)\in \ccalE\}$ of nodes that can communicate with $i$. In problem \eqref{original_optimization_problem1} agent $i$ has access to the local objective function $f_{i}(\tbx)$ and agents cooperate to minimize the global cost $\sum_{i=1}^nf_{i}(\tbx)$. This specification is more naturally formulated by defining variables $\bbx_i$ representing the local copies of the variable $\tbx$. In DADMM we further introduce the auxiliary variables $\bbz_{ij}$ associated with edge $(i,j)\in \ccalE$ and rewrite \eqref{original_optimization_problem1} as
\begin{align}\label{original_optimization_problem2}
 %  \{\bbx_i^*\}_{i=1}^n\ 
%   &:= \
%  \argmin_{\bbx} 
&\min_{\bbx,\bbz} \ \sum_{i=1}^{n}\ f_{i}(\bbx_{i}), \nonumber\\ 
   &\text{\ s.t.}  \ \bbx_{i}=\bbz_{ij},\ \bbx_{j}=\bbz_{ij}, 
                   \quad \text{for all\ } (i, j)\in\ccalE .
\end{align} 
The constraints $\bbx_{i}=\bbz_{ij}$ and $ \bbx_{j}=\bbz_{ij}$ enforce consensus among neighbors and, since the network is connected, they also enforce global consensus. 
Therefore, problems \eqref{original_optimization_problem1} and \eqref{original_optimization_problem2} are equivalent in the sense that for all $i$ and $j$ the optimal arguments of \eqref{original_optimization_problem2} satisfy $\bbx_i^*=\tbx^*$ for all $i\in\ccalV$ and $\bbz_{ij}=\tbx^*$, for all $(i,j)\in\ccalE$.

To specify the DADMM algorithm we write the constraints of \eqref{original_optimization_problem2} in matrix form. Begin by defining the block source matrix $\bbA_s\in \reals^{mp\times np}$ which contains $m\times n$ square blocks $(\bbA_s)_{e,i}\in \reals^{p\times p}$. The block $(\bbA_s)_{e,i}$ is not identically null if and only if the edge $e=(i,j)$ corresponds to $(i,j)\in \ccalE$ in which case $(\bbA_s)_{e,i}=\bbI_{p}$. Likewise, define the block destination matrix $\bbA_d\in \reals^{mp\times np} $ containing $m\times n$ square blocks $(\bbA_d)_{e,i}\in \reals^{p\times p}$. The square block $(\bbA_d)_{e,i}$ is equal to identity matrix  $\bbI_p$ if and only if $e=(j,i)$ corresponds to $(j,i)\in \ccalE$, otherwise the block is null. Further define the vector $\bbx:=[\bbx_1,\dots,\bbx_n]\in \reals^{np}$ as the vector that concatenates the local variables $\bbx_i$ and the vector $\bbz:=[\bbz_1,\dots,\bbz_m]\in \reals^{mp}$ that concatenates the auxiliary variables $\bbz_e=\bbz_{ij}$. If we further define the aggregate function $f:\reals^{np}\to\reals$ as $f(\bbx):= \sum_{i=1}^{n}f_i(\bbx_i)$ we can rewrite \eqref{original_optimization_problem2} as
\begin{align}\label{original_optimization_problem4}
   &\min_{\bbx,\bbz}f(\bbx), \qquad \st\ \bbA\bbx+\bbB\bbz=\bb0.
\end{align} 
where the matrices $\bbA\in \reals^{2mp\times np}$ and $\bbB\in \reals^{2mp\times mp}$ are defined as the stacks $\bbA:=[\bbA_s;\bbA_d]$ and $\bbB:=[-\bbI_{mp};-\bbI_{mp}]$.

Introduce now multipliers  $\bbalpha_e=\bbalpha_{ij}$ and $\bbbeta_e=\bbbeta_{ij}$ respectively associated with the constraints $\bbx_{i}=\bbz_{ij}$ and $\bbx_{j}=\bbz_{ij}$ [cf. \eqref{original_optimization_problem2}]. Concatenate the multipliers $\bbalpha_e$ in the vector $\bbalpha:=[\bbalpha_1,\dots,\bbalpha_m]$ and, likewise, stack the Lagrange multipliers $\bbbeta_e$ in the vector $\bbbeta:=[\bbbeta_1,\dots,\bbbeta_m]$. Further stacking the vectors $\bbalpha$ and $\bbbeta$ into $\bblambda:=[\bbalpha;\bbbeta]\in \reals^{2mp}$ leads to the Lagrange multiplier $\bblambda$ associated with the constraint $\bbA\bbx+\bbB\bbz=\bb0$ in \eqref{original_optimization_problem4}. The augmented Lagrangian of \eqref{original_optimization_problem4}, and \eqref{original_optimization_problem2}, which is the same optimization problem written with different notation, is now defined as
\begin{equation}\label{lagrangian}
\ccalL(\bbx,\bbz,\bblambda) := f(\bbx)+\bblambda^{T}\left(\bbA\bbx+\bbB\bbz\right)+ \frac{c}{2}\left\|\bbA\bbx+\bbB\bbz\right\|^2,
\end{equation}
where $c>0$ is an arbitrary positive constant.

The idea of the DADMM algorithm is to alternatively minimize the Lagrangian $\ccalL(\bbx,\bbz,\bblambda)$ with respect to the variable $\bbx$ and the auxiliary variable $\bbz$ and to follow these minimizations by a multiplier update collinear with the constraint violation. Specifically, consider a time index $k \in \naturals$ and define $\bbx_k$, $\bbz_k$, and $\bblambda_k$ as the primal and dual iterates at step $k$. The first step of DADMM is Lagrangian minimization with respect to $\bbx$ with $\bbz_k$ and $\bblam_k$ given,
\begin{equation}\label{ADMM_x_update}
   \bbx_{k+1}= \argmin_{\bbx}  f(\bbx)+\bblambda_k^{T}\left(\bbA\bbx+\bbB\bbz_k\right)
               + \frac{c}{2}\left\|\bbA\bbx+\bbB\bbz_k\right\|^2.
\end{equation}
The second step is minimization with respect to the auxiliary variable $\bbz$ with $\bblam_k$ given but using the updated variable $\bbx_{k+1}$,
\begin{equation}\label{ADMM_z_update}
\bbz_{k+1}\!=\! \argmin_{\bbz}  f(\bbx_{k+1})+\bblambda_k^{T}\left(\bbA\bbx_{k+1}\!+\!\bbB\bbz\right)+ \frac{c}{2}\left\|\bbA\bbx_{k+1}\!+\!\bbB\bbz\right\|^2\! .
\end{equation}
With the primal iterates $\bbx_{k+1}$ and $\bbz_{k+1}$ updated, the third and final step is to move the Lagrange multiplier $\bblambda_{k}$ in the direction of the constraint violation $\bbA\bbx_{k+1}+\bbB\bbz_{k+1}$,
\begin{equation}\label{ADMM_lambda_update}
\bblambda_{k+1}=\bblambda_{k}+c\left(\bbA\bbx_{k+1}+\bbB\bbz_{k+1}\right),
\end{equation}
where the constant $c$ is the same constant used in \eqref{lagrangian}. The DADMM algorithm follows from the observation that the updates in \eqref{ADMM_x_update}-\eqref{ADMM_lambda_update} can be implemented in a distributed manner, \cite{Schizas2008-1,BoydEtalADMM11,Shi2014-ADMM}. 

Notice that the update formulas in \eqref{ADMM_z_update} and \eqref{ADMM_lambda_update} have low computational cost. However, the update for the primal variable $\bbx$ in \eqref{ADMM_x_update} requires solution of an optimization problem. To avoid the cost of this minimization, the primal variable $\bbx_{k+1}$ can be updated inexactly. This idea leads to the DLM algorithm that approximates the objective function value $f(\bbx_{k+1})$ in \eqref{ADMM_x_update} through a linearization of the function $f$ in a neighborhood of the current variable $\bbx_{k}$. This approximation is defined as $f(\bbx_k)+\nabla f(\bbx_k)^T(\bbx-\bbx_k)+\rho\|\bbx-\bbx_k\|^2$, where $\rho$ is an arbitrary positive constant. Using this approximation, the update formula for the primal variable $\bbx$ in DLM replaces \eqref{ADMM_x_update} by 
\begin{align}\label{DLM_x_update}
\!\!\!\bbx_{k+1}= \argmin_{\bbx}\ & f(\bbx_k)+\nabla f(\bbx_k)^T(\bbx-\bbx_k)+\rho\|\bbx-\bbx_k\|^2
\nonumber \\
& \!\!\quad 
+\bblambda_k^{T}\left(\bbA\bbx+\bbB\bbz_k\right)+ \frac{c}{2}\left\|\bbA\bbx+\bbB\bbz_k\right\|^2.
\end{align}
Since the problem in \eqref{DLM_x_update} is quadratic, its solution is elementary. E.g., the first order optimality condition for \eqref{DLM_x_update} requires the variable $\bbx_{k+1}$ to satisfy
\begin{equation}\label{DLM_optimality_cond}
\nabla f(\bbx_{k})+\rho (\bbx_{k+1}-\bbx_{k})+\bbA^T\bblambda_{k}+c \bbA^{T} \left( \bbA\bbx_{k+1}+\bbB \bbz_{k}  \right)=\bb0.
\end{equation}
The expression in \eqref{DLM_optimality_cond} is a linear equation for $\bbx$ that an be solved by inversion of the positive definite matrix $\rho\bbI+c\bbA^T\bbA$. 

The sequence of variables $\bbx_k$ generated by DLM converges linearly to the optimal argument $\bbx^*$ \cite{cQingRibeiroADMM14}. This is the same rate of convergence of DADMM, but the linear rate coefficient of DLM is strictly smaller than the linear rate coefficient of DADMM \cite{DQMjournal}. In this paper we propose an alternative approximation that will be shown to achieve linear convergence with a coefficient that is asymptotically equivalent to the DADMM coefficient. This approximation utilizes second order information of the local functions $f_i(\bbx)$ and leads to the DQM algorithm that we introduce in the following section.

\section{DQM: Decentralized Quadratically Approximated ADMM}\label{sec:DQM}

We introduce a Decentralized Quadratic Approximation of ADMM (DQM) that uses a quadratic approximation of the primal function $f(\bbx)$ for updating the variable $\bbx$. To be more precise, instead of using $f(\bbx_{k+1})$ in the optimization problem \eqref{ADMM_x_update}, DQM executes the quadratic  approximation $f(\bbx_k)+\nabla f(\bbx_k)^T(\bbx-\bbx_k)+(1/2)(\bbx-\bbx_k)^T\bbH_k(\bbx-\bbx_k)$ for computing the new iterate $\bbx_{k+1}$, where $\bbH_{k}:=\nabla^2 f(\bbx_{k})$ is the Hessian of the primal function $f$ computed at $\bbx_k$ which is a block diagonal matrix. Therefore, the primal variable $\bbx$ is updated as
\begin{align}\label{DQM_x_update}
\bbx_{k+1}= \argmin_{\bbx} &\  f(\bbx_k) +\nabla f(\bbx_k)^T(\bbx-\bbx_k) \\
&\quad+\frac{1}{2}(\bbx-\bbx_k)^T\bbH_k(\bbx-\bbx_k)\nonumber \\
&\quad
+\bblambda_k^{T}\left(\bbA\bbx+\bbB\bbz_k\right)+ \frac{c}{2}\left\|\bbA\bbx+\bbB\bbz_k\right\|^2.\nonumber
\end{align}
The DQM updates for the variables $\bbz_k$ and $\bblambda_k$ are identical to the DADMM updates in \eqref{ADMM_z_update} and \eqref{ADMM_lambda_update}, respectively. 

Comparison of the updates in \eqref{DLM_x_update} and \eqref{DQM_x_update} shows that in DLM the quadratic term $\rho\|\bbx_{k+1}-\bbx_{k}\|^2$ is added to the first-order approximation of the primal objective function, while in DQM the second-order approximation of the primal objective function is used to reach a more accurate approximation for $f(\bbx_{k+1})$. First order optimality conditions of updates in \eqref{DQM_x_update}, \eqref{ADMM_z_update} and \eqref{ADMM_lambda_update} imply that the DQM iterates can be generated by solving the equations
\begin{align}\label{DQM_update}
\nabla f(\bbx_{k})+ & \bbH_{k}(\bbx_{k+1}-\bbx_{k})+\bbA^T\!\bblambda_{k}\!+\!c \bbA^{T}\! \left( \bbA\bbx_{k+1}\!+\!\bbB \bbz_{k}  \right)=\bb0,\nonumber\\
&\bbB^T\bblambda_{k} +c \bbB^T\left( \bbA\bbx_{k+1}+\bbB\bbz_{k+1} \right) =\bb0,\nonumber\\
&\bblambda_{k+1}-\bblambda_{k}-c\left( \bbA\bbx_{k+1}+\bbB\bbz_{k+1} \right) =\bb0.
\end{align}
Notice that the first equation in \eqref{DQM_update} can be solved by inverting $\bbH_k+c\bbA^T\bbA$. To guarantee that the system of equations in \eqref{DQM_update} can be solved in a distributed manner, we assume a specific structure for the initial vectors $\bblambda_{0}=[\bbalpha_0;\bbbeta_0]$, $\bbx_0$, and $\bbz_0$ as mentioned in Assumption \ref{initial_val_assum}. Before introducing this assumption we define the oriented incidence matrix as $\bbE_o:=\bbA_s-\bbA_d$ and the unoriented incidence matrix as $\bbE_u:=\bbA_s+\bbA_d$.

%%%%%%%%%%%%%%%%%%%%%%%%%%%%%%%%%%%%%%%%%%%%
%%%%%%%%%%%%%%%%%%%%%%%%%%%%%%%%%%%%%%%%%%%%
%%%%%%     A   S   S   U   M   P   T   I   O   N     %%%%%%%%%%%%%%%%%%%
%%%%%%%%%%%%%%%%%%%%%%%%%%%%%%%%%%%%%%%%%%%%
%%%%%%%%%%%%%%%%%%%%%%%%%%%%%%%%%%%%%%%%%%%%
\begin{assumption}\label{initial_val_assum}

The initial Lagrange multipliers $\bbalpha_0$ and $\bbbeta_0$, and the initial variables $\bbx_0$ and $\bbz_0$ are chosen such that
\begin{enumerate}[a)]
\item $\bbalpha_0=-\bbbeta_0$,
\item $\bbE_u\bbx_0=2\bbz_0$,
\item $\bbalpha_0$ lies in the column space of $\bbE_o$.
\end{enumerate}
%The initial Lagrange multipliers satisfy $\bbalpha_0=-\bbbeta_0$. Further, the initial primal variable $\bbx_0$ and auxiliary variable $\bbz_0$ are chosen such that $\bbE_u\bbx_0=2\bbz_0$. The initial Lagrangian vector $\bbalpha_0$ lies in the column vector of the oriented incidence matrix $\bbE_o$, i.e. there exits a vector $\bbr_0$ such that $\bbalpha_0=\bbE_o\bbr_0$.
\end{assumption}

%%%%%%%%%%%%%%%%%%%%%%%%%%%%%%%%%%%%%%%%%
%%%%%%%%%%%%%%%%%%%%%%%%%%%%%%%%%%%%%%%%%
%%%%%%%%%%%%%%%%%%%%%%%%%%%%%%%%%%%%%%%%%
%%%   M   A   I   N       M   A   T   T   E   R   %%%%%%%%%%%%%%%%%%%
%%%%%%%%%%%%%%%%%%%%%%%%%%%%%%%%%%%%%%%%%
%%%%%%%%%%%%%%%%%%%%%%%%%%%%%%%%%%%%%%%%%

Assumption \ref{initial_val_assum} enforces some initial conditions on the Lagrange multipliers $\bbalpha_0$ and $\bbbeta_0$, and the initial variables $\bbx_0$ and $\bbz_0$. These conditions are satisfied by setting $\bbalpha_0=\bbbeta_0=\bb0$ and $\bbx_0=\bbz_0=\bb0$. 
%In the following lemma we show that satisfaction of these conditions for the initial values guarantees the same conditions for all steps $k>0$. 
%
%%\blue{bayad ezafe konam inja, ye lemma ke shamele 3 bakhsh hast}\\
%%\red{ Should I mention Lemma 3 and 4 right after Assumption 1 to explain the relation between $\bbalpha $ and $\bbbeta$, and the relation between $\bbx$ and $\bbz$?}
%
%
%%%%%%%%%%%%%%%%%%%%%%%%%%%%%%%%%%%%%%%%%
%%%%%%%%%%%%%%%%%%%%%%%%%%%%%%%%%%%%%%%%%
%%%%%%%    L   E   M   M   A    %%%%%%%%%%%%%%%%%%%%%%%%
%%%%%%%%%%%%%%%%%%%%%%%%%%%%%%%%%%%%%%%%%
%%%%%%%%%%%%%%%%%%%%%%%%%%%%%%%%%%%%%%%%%
%\begin{lemma}\label{lag_var_lemma}
%Consider the DQM algorithm as defined in \eqref{DQM_x_update} and \eqref{DQM_update}. If Assumption \ref{initial_val_assum} holds true, then for $k\geq 0$ the Lagrange multipliers $\bbalpha_k$ and $\bbbeta_k$, and the variables $\bbx_k$ and $\bbz_k$ satisfy 
%\begin{enumerate}[a)]
%\item $\bbalpha_k=-\bbbeta_k$,
%\item $\bbE_u\bbx_k=2\bbz_k$,
%\item $\bbalpha_k$ lies in the column space of $\bbE_o$.
%\end{enumerate}
%%%
%\end{lemma}
%
%%%%%%%%%%%%%%%%%%%%%%%%%%%%%%%%%%%%%%%%%%
%%%%%%%%%%%%%%%%%%%%%%%%%%%%%%%%%%%%%%%%%%
%%%%%%%%%%%%%%%%%%%%%%%%%%%%%%%%%%%%%%%%%%
%%%%   M   A   I   N       M   A   T   T   E   R   %%%%%%%%%%%%%%%%%%%
%%%%%%%%%%%%%%%%%%%%%%%%%%%%%%%%%%%%%%%%%%
%%%%%%%%%%%%%%%%%%%%%%%%%%%%%%%%%%%%%%%%%%
%
%Lemma \ref{lag_var_lemma} states that if the conditions in Assumption \ref{initial_val_assum} are satisfied for $k=0$, the same relationships are also valid for all steps $k\geq0$. 
Considering the initial conditions in Assumption  \ref{initial_val_assum} we propose a simpler update rule for generating the iterates $\bbx_k$ instead of solving equations in \eqref{DQM_update}.

%%%%%%%%%%%%%%%%%%%%%%%%%%%%%%%%%%%%%%%%
%%%%%%%%%%%%%%%%%%%%%%%%%%%%%%%%%%%%%%%%
%%%%%%    P   R   O  P  O  S   I    T   I   O   N    %%%%%%%%%%%%%%%
%%%%%%%%%%%%%%%%%%%%%%%%%%%%%%%%%%%%%%%%
%%%%%%%%%%%%%%%%%%%%%%%%%%%%%%%%%%%%%%%%
\begin{proposition}\label{update_system_prop}
Consider the system of equations for the DQM algorithm in \eqref{DQM_update} and define the sequence $\bbphi_{k}:=\bbE_{o}^T\bbalpha_k$. Further, define the unoriented Laplacian as $\bbL_u:=(1/2)\bbE_u^T\bbE_u$, the oriented Laplacian as $\bbL_o:=(1/2)\bbE_o^T\bbE_o$, and the Degree matrix as $\bbD:=(\bbL_u+\bbL_o)/2$. If Assumption \ref{initial_val_assum} holds true, the DQM variables $\bbx_k$ can be generated as
\begin{equation}\label{x_update_formula}
\!\!\bbx_{k+1}=(2c\bbD+\bbH_k)^{-1}\left[ (c\bbL_u+\bbH_k)\bbx_k-\nabla f(\bbx_{k})-\bbphi_k \right],
\end{equation}
\begin{equation}\label{phi_update_formula}
\bbphi_{k+1} =\bbphi_k +c \bbL_o\bbx_{k+1}.
\end{equation}
\end{proposition}

%%%%%%%%%%%%%%%%%%%%%%%%%%%%%%%%%%%%%%%%%
%%%%%%%%%%%%%%%%%%%%%%%%%%%%%%%%%%%%%%%%%
%%%%%%%%%%%%%%%%%%%%%%%%%%%%%%%%%%%%%%%%%
%%%   M   A   I   N       M   A   T   T   E   R   %%%%%%%%%%%%%%%%%%%
%%%%%%%%%%%%%%%%%%%%%%%%%%%%%%%%%%%%%%%%%
%%%%%%%%%%%%%%%%%%%%%%%%%%%%%%%%%%%%%%%%%
 
 Proposition \ref{update_system_prop} states that by introducing the new variables $\bbphi_{k}$, the update formulas for the DQM iterates can be computed using the primal objective function Hessian $\bbH_{k}$, the degree matrix $\bbD$, and the oriented and unoriented Laplacians $\bbL_o$ and $\bbL_u$. This observation guarantees that the updates in \eqref{x_update_formula} and \eqref{phi_update_formula} are implementable in a distributed manner, since all of these matrices can be computed using local and neighboring information of the nodes. To be more precise, the matrix $2c\bbD+\bbH_k$ is block diagonal and its $i$-th diagonal block is given by $2cd_{i}\bbI+\nabla^2f_{i}(\bbx_{i})$ which is locally available at node $i$. Likewise, the inverse matrix $(2c\bbD+\bbH_k)^{-1}$ is block diagonal and locally computable since the $i$-th diagonal block is $(2cd_{i}\bbI+\nabla^2f_{i}(\bbx_{i}))^{-1}$. Computations of the products $\bbL_u\bbx_{k}$ and $\bbL_o\bbx_{k+1}$ can be implemented in a decentralized manner as well, since the Laplacian matrices $\bbL_u$ and $\bbL_o$ are block neighbor sparse. Note that a matrix is block neighbor sparse when its ${ij}$-th block is not null if and only if nodes $i$ and $j$ are neighbors or $j=i$. Therefore, nodes can compute $\bbL_u\bbx_{k}$ and $\bbL_o\bbx_{k+1}$ by exchanging information with their neighbors. By defining the components of $\bbphi_k$ as $\bbphi_k:=[\bbphi_{1,k},\dots,\bbphi_{n,k}]$, the update in \eqref{x_update_formula} can be implemented locally as 
 \begin{align}\label{x_local_update_formula}
\bbx_{i,k+1}=\ & \left(2cd_{i}\bbI+\nabla^2f_{i}(\bbx_{i,k})\right)^{-1}\Big[ cd_{i}\bbx_{i,k}+c\sum_{j\in \ccalN_i}\bbx_{j,k}
\nonumber \\
&\qquad +\nabla^2f_{i}(\bbx_{i,k})\bbx_{i,k}-\nabla f_i(\bbx_{i,k})-\bbphi_{i,k} \Big],
\end{align}
where $\bbx_{i,k}$ and $\bbphi_{i,k}$ are the iterates of node $i$ at step $k$.
Notice that the definition $\bbL_u:=(1/2)\bbE_u^T\bbE_u=(1/2)(\bbA_s+\bbA_d)^T(\bbA_s+\bbA_d)$ is used to simplify the $i$-th component of $c\bbL_u\bbx_{k}$ as $cd_{i}\bbx_{i,k}+c\sum_{j\in \ccalN_i}\bbx_{j,k}$. Further, considering the definition $\bbL_o=(1/2)\bbE_o^T\bbE_o=(1/2)(\bbA_s-\bbA_d)^T(\bbA_s-\bbA_d)$, the $i$-th component of $c\bbL_o\bbx_{k+1}$ can be simplified as $c\sum_{j\in \ccalN_i}(\bbx_{i,k}-\bbx_{j,k})$. Therefore, the update formula in \eqref{phi_update_formula} for node $i $ is given by
\begin{equation}\label{phi_local_update_formula}
\bbphi_{i,k+1} =\bbphi_{i,k} +c \sum_{j\in \ccalN_i} \left(\bbx_{i,k+1} -\bbx_{j,k+1}\right).
\end{equation}
The proposed DQM algorithm is summarized in Algorithm \ref{algo_DQM}. 
%The initial value for the local iterate $\bbx_{i,0}$ can be any arbitrary vector in $\reals^p$. The initial vector $\bbphi_{i,0}$ should be in the column space of $\bbE_o^T$. To guarantee satisfaction of this condition, the initial vector is set as $\bbphi_{i,0}=\bb0$. 
At each iteration $k$, the primal and dual updates in \eqref{x_local_update_formula} and \eqref{phi_local_update_formula} are implemented in Steps 2 and 4, respectively. Nodes exchange their local variables $\bbx_{i,k}$ with their neighbors $j\in \ccalN_i$ in Step 3, since this information is required for Steps 2 and 4.
%%%%%%%%%%%%%%%%%%%%%%%%%%%%%%%%%%%%%%%%%%%%%%%%%%%%%%%%%%%%%%%%%%
%%%   A   L   G   O   R   I   T   H   M   %%%%%%%%%%%%%%%%%%%%%%%%
%%%%%%%%%%%%%%%%%%%%%%%%%%%%%%%%%%%%%%%%%%%%%%%%%%%%%%%%%%%%%%%%%%
%
\begin{algorithm}[t]{\small
\caption{DQM method at node $i$}\label{algo_DQM} 
\begin{algorithmic}[1] {
\REQUIRE  Initial local iterates $\bbx_{i,0}$ and $\bbphi_{i,0}$. 
%\STATE $\bbB$ matrix blocks: $\bbB_{ii}=(1-w_{ii})\bbI$ and $\bbB_{ij}=w_{ij}\bbI$
\FOR {$k=0,1,2,\ldots$}

   \STATE Update the local iterate $\bbx_{i,k+1}$
    \begin{align}
\bbx_{i,k+1}=\ & \left(2cd_{i}\bbI+\nabla^2f_{i}(\bbx_{i,k})\right)^{-1}\bigg[ cd_{i}\bbx_{i,k}+c\sum_{j\in \ccalN_i}\bbx_{j,k}
\nonumber \\
&\qquad+\nabla^2f_{i}(\bbx_{i,k})\bbx_{i,k}-\nabla f_i(\bbx_{i,k})-\bbphi_{i,k} \bigg].\nonumber
\end{align}

   \STATE Exchange iterates $\bbx_{i,k+1}$ with neighbors $\displaystyle{j\in \mathcal{N}_i}$.
   
   \STATE Update local dual variable $\bbphi_{k+1} $ as 
         \begin{equation}
\bbphi_{i,k+1} =\bbphi_{i,k} +c \sum_{j\in \ccalN_i} \left(\bbx_{i,k+1} -\bbx_{j,k+1}\right). \nonumber
\vspace{-2mm}
\end{equation}
\ENDFOR}
\end{algorithmic}}\end{algorithm}

%%%%%%%%%%%%%%%%%%%%%%%%%%%%%%%%%%%%%%%%
%%%%%%%%%%%%%%%%%%%%%%%%%%%%%%%%%%%%%%%%
%%%%%%     R  E   M  A  R  K    %%%%%%%%%%%%%%%%%%%%%%%
%%%%%%%%%%%%%%%%%%%%%%%%%%%%%%%%%%%%%%%%
%%%%%%%%%%%%%%%%%%%%%%%%%%%%%%%%%%%%%%%%
%
%\begin{remark}
%Implementation of the DQM algorithm requires exchanging local iterates $\bbx_{i,k}$ with neighboring nodes at each iteration which is identical to the communication cost per iteration of DLM and DADMM. Though the communication costs per iteration of these three algorithms are equal, the computation complexity of DLM and DQM algorithms are significantly lower relative to DADMM. DADMM requires minimizing a quadratic convex program for each node at each iteration, while the computation complexities of DLM and DQM are in the order of $\ccalO(p)$ and $\ccalO(p^2)$, respectively. DQM computation complexity per iteration is higher than $DLM$, although, the runtime until convergence for DQM is less than the one for DLM as DQM requires less iterations to converge relative to DLM. These results are shown both theoretically and numerically -- see Sections \ref{sec:Analysis} and \ref{sec:simulations}.
%\end{remark}

%!TEX root = root.tex

\section{Convergence Analysis} \label{sec:Analysis}

In this section we show that the sequence of iterates $\bbx_k$ generated by the DQM method converges linearly to the optimal argument $\bbx^*$. In addition, we compare the linear convergence coefficients of the DQM and DADMM methods. To provide these results first we make the following assumptions.

%%%%%%%%%%%%%%%%%%%%%%%%%%%%%%%%%%%%%
%%%%%%%%%%%%%%%%%%%%%%%%%%%%%%%%%%%%%
%%%   A   S   S   U   M   P   T   I   O   N    %%%%%%%%%%%%%%%%
%%%%%%%%%%%%%%%%%%%%%%%%%%%%%%%%%%%%%
%%%%%%%%%%%%%%%%%%%%%%%%%%%%%%%%%%%%%

\begin{assumption}\label{convexity_assumption} 
The local objective functions $f_i(\bbx)$ are twice differentiable and the eigenvalues of the local Hessians $ \nabla^2 f_i(\bbx)$ are bounded with positive constants $0<m\leq M<\infty$. I.e., 
\begin{equation}\label{local_hessian_eigenvlaue_bounds}
m\bbI\preceq \nabla^2 f_i(\bbx)\preceq M\bbI.
\end{equation}
\end{assumption}

%%%%%%%%%%%%%%%%%%%%%%%%%%%%%%%%%%%%%
%%%%%%%%%%%%%%%%%%%%%%%%%%%%%%%%%%%%%
%%%   A   S   S   U   M   P   T   I   O   N    %%%%%%%%%%%%%%%%
%%%%%%%%%%%%%%%%%%%%%%%%%%%%%%%%%%%%%
%%%%%%%%%%%%%%%%%%%%%%%%%%%%%%%%%%%%%

\begin{assumption}\label{Lipschitz_assumption} The local Hessians $\nabla^2 f_i(\bbx)$ are Lipschitz continuous with parameter $L$. I.e., for all $\bbx, \hbx \in \reals^p$, it holds
\begin{equation}
   \left\| \nabla^2 f_i(\bbx)-\nabla^2 f_i(\hbx) \right\| \ \leq\  L\ \| \bbx- \hbx \|.
\end{equation}

\end{assumption}

%%%%%%%%%%%%%%%%%%%%%%%%%%%%%%%%%%%%%%%%%
%%%%%%%%%%%%%%%%%%%%%%%%%%%%%%%%%%%%%%%%%
%%%   M   A   I   N       M   A   T   T   E   R   %%%%%%%%%%%%%%%%%%%
%%%%%%%%%%%%%%%%%%%%%%%%%%%%%%%%%%%%%%%%%
%%%%%%%%%%%%%%%%%%%%%%%%%%%%%%%%%%%%%%%%%

{The lower bound $m$ for the eigenvalues of the local Hessians $\nabla^2 f_i(\bbx)$ implies that the local objective functions $f_{i}(\bbx)$ are strongly convex with parameter $m$. The upper bound $M$ for the eigenvalues of the local Hessians $\nabla^2 f_i(\bbx)$ is similar to the condition that the local gradients $\nabla f_i(\bbx)$ are Lipschitz continuous with parameter $M$. The Lipschitz continuity of the local Hessians imposed by Assumption \ref{Lipschitz_assumption} is typical of second order methods.

 Define $\gamma_u$ and $\Gamma_u$ as the minimum and maximum singular values of the unoriented incidence matrix $\bbE_u$, respectively. Further, define $\gamma_o$ as the smallest non-zero singular value of the oriented incidence matrix $\bbE_o$. These parameters capture connectivity of the network. Denote the unique solution of \eqref{original_optimization_problem2} as $(\bbx^*,\bbz^*)$. Notice that the uniqueness is implied by the strong convexity assumption. Further, define $\bbalpha^*$ as the unique optimal multiplier that lies in the column space of $\bbE_o$ -- see Lemma 1 of \cite{ling2014dlm} for the uniqueness of such an optimal dual variable $\bbalpha^*$. We define the energy function $V: \reals^{mp\times mp}\to\reals$ as introduced in \cite{Shi2014-ADMM} for DADMM,
\begin{equation}
V(\bbz,\bbalpha):=c\|\bbz-\bbz^{*}\|^2+(1/c)\|\bbalpha-\bbalpha^*\|^2.
\end{equation}
The energy function $V(\bbz,\bbalpha)$ captures the distances of the auxiliary variable $\bbz_{k}$ and the dual variable $\bbalpha_{k}$ with their optimal arguments $\bbz^{*} $ and $\bbalpha^*$, respectively. Therefore, convergence rate of the energy function is a valid tool for comparing performances of DQM and DADMM. To simplify the notation of the energy function $V(\bbz,\bbalpha)$, define $\bbu\in \reals^{2mp}$ and $\bbC\in \reals^{2mp\times 2mp}$ as
%%%
\begin{equation}\label{C_u_definitions}
\bbu:= \left[ \begin{array}{cc}
\bbz  \\
\bbalpha  \end{array} \right], \quad 
\bbC:= \left[ \begin{array}{cc}
c\bbI & \bb0  \\
\bb0 & (1/c)\bbI \end{array} \right].
\end{equation}
The energy function $V(\bbz,\bbalpha)$ is equal to weighted squared norm $\|\bbu-\bbu^*\|_{\bbC}^2$ where $\bbu^*:=[\bbz^*;\bbalpha^*]$. Our goal is to show that the sequence $\|\bbu_{k}-\bbu^*\|_\bbC^2$ converges linearly to null.

%%%%%%%%%%%%%%%%%%%%%%%%%%%%%%%%%%%%%%%%%%
%%%%%%%%%%%%%%%%%%%%%%%%%%%%%%%%%%%%%%%%%%
%%%%%%%%%%%%%%%%%%%%%%%%%%%%%%%%%%%%%%%%%%
%%%%%%%%%%     T  H  E  O  R  E  M      %%%%%%%%%%%%%%%%%%%
%%%%%%%%%%%%%%%%%%%%%%%%%%%%%%%%%%%%%%%%%%
%%%%%%%%%%%%%%%%%%%%%%%%%%%%%%%%%%%%%%%%%%
%%%%%%%%%%%%%%%%%%%%%%%%%%%%%%%%%%%%%%%%%%
%%
\begin{thm}\label{DQM_convergence}
Consider the DQM method as introduced in \eqref{DQM_x_update}-\eqref{phi_local_update_formula}. Define the sequence of non-negative variables $\zeta_k$ as
\begin{equation}\label{zeta_DQM}
\zeta_k:=\min\left\{\frac{L}{2}\|\bbx_{k+1}-\bbx_{k} \|, 2M\right\}.
\end{equation}
Assume that the constant $c$ is chosen such that $c> {\zeta_k^2}/({m\gamma_u^2})$.
Moreover, consider $\mu,\mu'>1$ as arbitrary constants and $\eta$ as a positive constant chosen from the interval $({\zeta_k}/{m},{c\gamma_u^2}/{\zeta_k})$.
If Assumptions \ref{initial_val_assum}-\ref{Lipschitz_assumption} hold true, then the sequence $\|\bbu_{k}-\bbu^*\|_{\bbC}^2$ generated by DQM satisfies 
%%%
\begin{equation}\label{DQM_linear_claim}
\|\bbu_{k+1}-\bbu^*\|_{\bbC}^2\  \leq\ \frac{1}{1+\delta_k}  \|\bbu_{k}-\bbu^*\|_{\bbC}^2\ ,
\end{equation}
where the sequence of positive scalars $\delta_k$ is given by
\begin{equation}\label{DLM_DQM_delta}
\delta_k\!=\!\min\!
\Bigg\{\!
\frac{c-\eta\zeta_k\gamma_u^{-2} }{\frac{4\mu'\mu\zeta_k^2}{c(\mu'\!-\!1)(\mu-1)}\gamma_u^{-2}\gamma_o^{-2}\!+\!\frac{\mu'\mu}{(\mu-1)}\Gamma_u^2\gamma_o^{-2}}
,
\frac{  m-{\zeta_k}/{\eta} }{\frac{ c}{4}\Gamma_u^2
\!+\!\frac{\mu }{c}M^2\gamma_o^{-2}}
\!\Bigg\}\!.
\end{equation}
%%%
\end{thm}

%%%%%%%%%%%%%%%%%%%%%%%%%%%%%%%%%%%%%%%%%
%%%%%%%%%%%%%%%%%%%%%%%%%%%%%%%%%%%%%%%%%
%%%%%%%%%%%%%%%%%%%%%%%%%%%%%%%%%%%%%%%%%
%%%   M   A   I   N       M   A   T   T   E   R   %%%%%%%%%%%%%%%%%%%
%%%%%%%%%%%%%%%%%%%%%%%%%%%%%%%%%%%%%%%%%
%%%%%%%%%%%%%%%%%%%%%%%%%%%%%%%%%%%%%%%%%

Notice that $\delta_k$ is a decreasing function of $\zeta_k$ and observe that $\zeta_k$ is bounded above by $2M$. Therefore, if we substitute $\zeta_k$ by $2M$ in \eqref{DLM_DQM_delta}, the inequality in \eqref{DQM_linear_claim} is still valid. This substitution implies that the sequence $\|\bbu_{k}-\bbu^*\|_{\bbC}^2$ converges linearly to null. %Note that in order to guarantee that $\delta_k>0$ for all $k\geq 0 $, the constant $\eta$ is chosen from the interval $({\zeta_k}/{m},{c\gamma_u^2}/{\zeta_k})$. This interval is non-empty if the constant $c$ is chosen as $c> {\zeta_k^2}/({m\gamma_u^2})$. 
As a result of this convergence we obtain that $\bbu_{k}$ approaches the optimal argument $\bbu^*$. Therefore, the sequence of primal iterates $\bbx_k$ converges to the optimal argument $\bbx^*$. This result is formalized in the following corollary.

%%%%%%%%%%%%%%%%%%%%%%%%%%%%%%%%%%%%%%%%%%
%%%%%%%%%%%%%%%%%%%%%%%%%%%%%%%%%%%%%%%%%%
%%%%%%%%%%%%%%%%%%%%%%%%%%%%%%%%%%%%%%%%%%
%%%%%%%%%%     C  O   R   O    L     L     A     R    Y      %%%%%%%%%%%%
%%%%%%%%%%%%%%%%%%%%%%%%%%%%%%%%%%%%%%%%%%
%%%%%%%%%%%%%%%%%%%%%%%%%%%%%%%%%%%%%%%%%%
%%%%%%%%%%%%%%%%%%%%%%%%%%%%%%%%%%%%%%%%%%
%%
\begin{corollary}\label{convergence_of_primal}
Under the assumptions in Theorem \ref{DQM_convergence}, the sequence of squared norms $\|\bbx_k-\bbx^*\|^2$ generated by the DQM algorithm converges R-linearly to null, i.e.,
%%%
\begin{equation}\label{r_linear_cliam}
\|\bbx_k-\bbx^*\|^2\leq \frac{4}{c\gamma_u^2}\|\bbu_k-\bbu^*\|_\bbC^2.
\end{equation}
\end{corollary}

%%%%%%%%%%%%%%%%%%%%%%%%%%%%%%%%%%%%%%%%%
%%%%%%%%%%%%%%%%%%%%%%%%%%%%%%%%%%%%%%%%%
%%%%%%%%%%%%%%%%%%%%%%%%%%%%%%%%%%%%%%%%%
%%%   M   A   I   N       M   A   T   T   E   R   %%%%%%%%%%%%%%%%%%%
%%%%%%%%%%%%%%%%%%%%%%%%%%%%%%%%%%%%%%%%%
%%%%%%%%%%%%%%%%%%%%%%%%%%%%%%%%%%%%%%%%%

Corollary \ref{convergence_of_primal} states that the sequence $\bbx_k$ converges to the optimal argument $\bbx^*$. Hence, we obtain that the sequence $\|\bbx_{k+1}-\bbx_k\|$ approaches null. This observation implies that the sequence of scalars $\zeta_k$ converges to 0 as time passes, since $\zeta_k$ is bounded above by $(L/2)\|\bbx_{k+1}-\bbx_k\|$. By considering $\lim_{k\to \infty }\zeta_k=0$ and making $\mu'$ arbitrary close to $1$, we obtain 
\begin{align}\label{DQM_delta}
\lim_{k\to \infty} \delta_k=\min
&\Bigg\{
\frac{ (\mu-1)\gamma_o^2}{ {\mu\Gamma_u^2}}
,\
\frac{  m}{\frac{ c}{4}\Gamma_u^2
+\frac{\mu}{c}M^2\gamma_o^{-2}}
\Bigg\}.
\end{align}
Notice that the coefficient $\delta_k$ in \eqref{DQM_delta} is identical to the coefficient of linear convergence for the DADMM algorithm \cite{Shi2014-ADMM}. This observation implies that as time passes the coefficient of linear convergence for DQM approaches the one for DADMM.

%\input{Convergence.tex}
%\input{Implementation.tex}
%!TEX root = root.tex

%%%%%%%%%%%%%%%%%%%%%%%%%%%%%%%%%%%%%%%%%%%%%%%%%%%%%%%%%%%%%%%%%%%%
%   S   E   C   T   I   O   N   %%%%%%%%%%%%%%%%%%%%%%%%%%%%%%%%%%%%
%%%%%%%%%%%%%%%%%%%%%%%%%%%%%%%%%%%%%%%%%%%%%%%%%%%%%%%%%%%%%%%%%%%%
%
\section{Numerical analysis}\label{sec:simulations}

 In this section we compare performances of DLM, DQM and DADMM in solving a logistic regression problem. We assume that each node in the network has access to $q$ training points. Therefore, the total number of training points is $nq$. Each of the training samples $\{\bbs_{il}, y_{il}\}_{l=1}^q$ at node $i$ contains a feature vector $\bbs_{il}\in \reals^p$ with  class $y_{il}\in \{-1,1\}$. It follows from the logistic regression model that the maximum log-likelihood estimate of the classifier $\tbx$ given the training samples $(\bbs_{il},y_{il})$ for $l=1,\ldots,q$ and $i=1,\ldots, n$ is 
\begin{align}\label{eqn_logistic_regrssion_max_likelihood}
   \tbx^* \ :=\  \argmin_{\tbx\in \reals^p }  \    \sum_{i=1}^n \sum_{l=1}^{q} 
                                  \log \Big[1+\exp(-y_{il}\bbs_{il}^T\tbx)\Big].
\end{align}
%
% The objective function in \eqref{eqn_logistic_regrssion_max_likelihood} is twice differentiable and strongly convex if every entry of $\bbs_il$ is finite
The optimization problem in \eqref{eqn_logistic_regrssion_max_likelihood} can be written in the form of \eqref{original_optimization_problem1} by defining the local objective functions $f_i$ as 
\begin{equation}
   f_i(\tbx) :=     \sum_{l=1}^{q} \log \Big[1+\exp(-y_{il}\bbs_{il}^T\tbx)\Big].
\end{equation}
%
%Therefore, we can solve \eqref{eqn_logistic_regrssion_max_likelihood} in a distributed manner. %We define the optimal argument for the decentralized optimization as $\bbx^*=[\tbx^*;\dots;\tbx^*]$.

%%%%%%%%%%%%%%%%%%%%%%%%%%%%%%%%%%%%%%%%%%%%%%%%%%%%%%%%%%%%%%%%%%%%
%   F   I   G   U   R   E   %%%%%%%%%%%%%%%%%%%%%%%%%%%%%%%%%%%%%%%%
%%%%%%%%%%%%%%%%%%%%%%%%%%%%%%%%%%%%%%%%%%%%%%%%%%%%%%%%%%%%%%%%%%%%
%
\begin{figure}[t]
\centering
\includegraphics[width=\linewidth,height=0.5\linewidth]{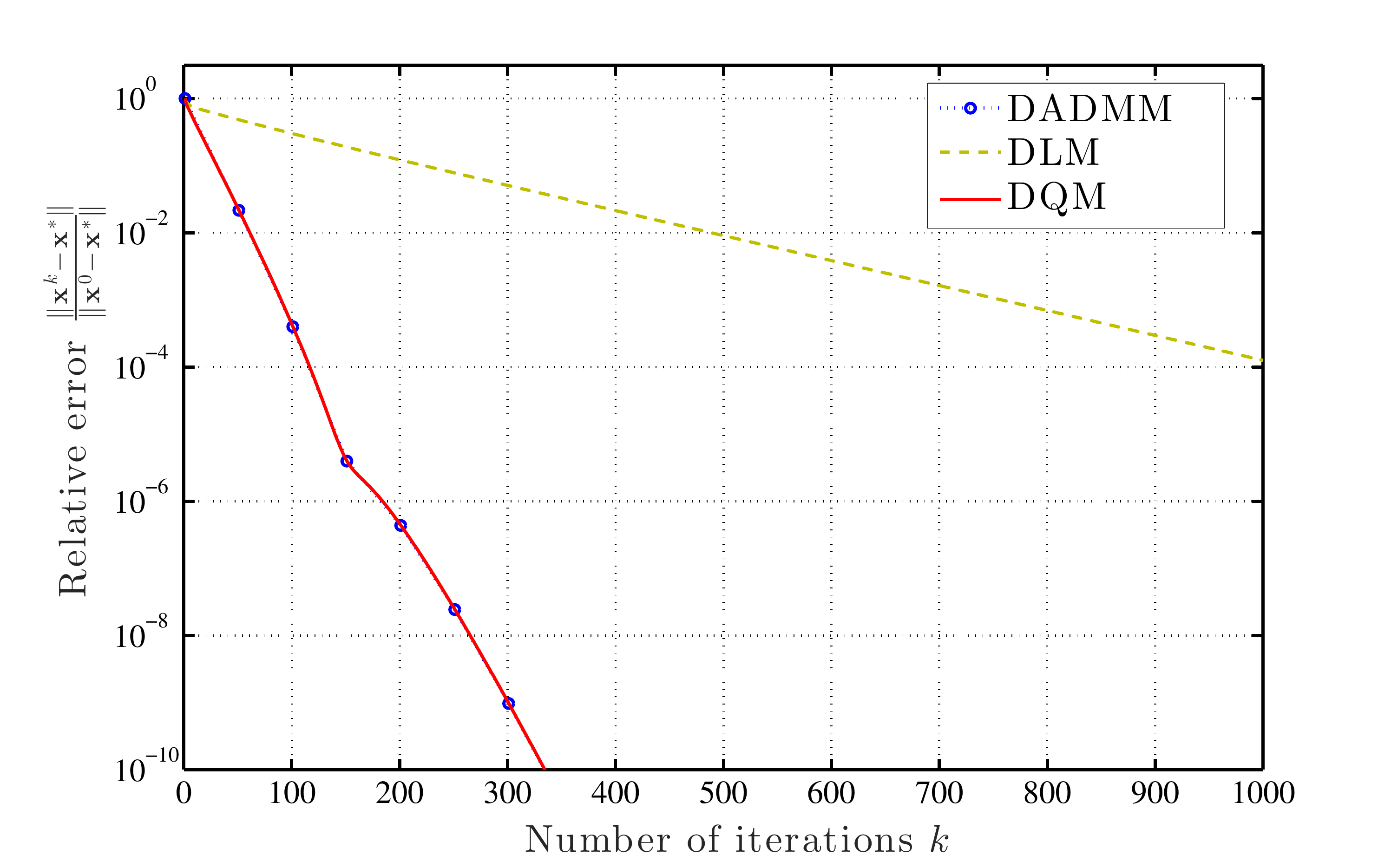}
\vspace{1mm}
\caption{{Relative error ${\|\bbx^k-\bbx^*\|}/{\|\bbx^0-\bbx^*\|}$ of DADMM, DQM, and DLM versus number of iterations. The convergence path of DQM is almost identical to the one for DADMM, while DQM has less computation complexity. Further, DQM outperforms DLM in convergence speed by orders of magnitude, though at the cost of higher computation complexity. }}
\label{fig2}
\vspace{2mm}
\end{figure}

We compare convergence paths of the DLM, DQM and DADMM algorithms for solving the logistic regression problem. We assume the network contains $n=10$ nodes and the edges between nodes are generated randomly with probability  $P_{\mathrm{c}}=0.4$. Each agent holds $q=5$ samples and the dimension of feature vectors is $p=3$. The reference (ground true) logistic classifier $\tbx^*$ is pre-computed with a centralized method. Notice that the parameter $c$ for the three methods is optimized by $c_\mathrm{ADMM}=0.7$, $c_\mathrm{DLM}=5.5$, and $c_\mathrm{DQM}=0.7$. Figure \ref{fig2} illustrates the relative errors ${\|\bbx^k-\bbx^*\|}/{\|\bbx^0-\bbx^*\|}$ of DLM, DQM, and DADMM versus the number of iterations. The convergence path of DQM is almost identical to the convergence path of DADMM and they both converge to the optimal argument faster than DLM. The relative errors ${\|\bbx^k-\bbx^*\|}/{\|\bbx^0-\bbx^*\|}$ for DQM and DADMM after $k=300$ iterations are below $10^{-9}$, while for DLM the relative error after the same number of iterations is $5\times 10^{-2}$. Conversely, achieving error  ${\|\bbx^k-\bbx^*\|}/{\|\bbx^0-\bbx^*\|}=10^{-3}$ for DQM and DADMM requires $91$ iterations, while DLM requires $758$ iterations. Observe that the convergence paths of DQM and DADMM are almost identical, while the computation complexity of DQM is lower than DADMM.

\newpage
\bibliographystyle{IEEEtran}
  \bibliography{bmc_article}
   \end{document}